\documentstyle{amsppt}
\NoBlackBoxes
\def\n{\noindent}

 2
\font\tenbb=msbm10
\font\sevenbb=msbm7
\font\fivebb=msbm5
\newfam\bbfam
\textfont\bbfam=\tenbb
\scriptfont\bbfam=\sevenbb
\scriptscriptfont\bbfam=\fivebb
\def\bb{\fam\bbfam\tenbb}
\font\tenll=lasy10
\newfam\llfam
\textfont\llfam=\tenll

\def\dc{\mathop{\nabla}}

\def\dczer{\dc\limits^o}

\def\dca{\dc\limits^{\alpha}}

\def\tr{\mathop{\triangle}}
\def\tra{\tr\limits^{\alpha}}

\leftheadtext{ Guy BURDET}
\rightheadtext{ c$\cup$p-invariant differential operators}

\topmatter
\title  INVARIANTS IN  c $\cup$ p-GEOMETRY  \endtitle
\author  Guy BURDET \endauthor
\address  Centre de Physique Th\'eorique, CNRS-Luminy, Case 907,
F-13288 Marseille, France\endaddress  
%Research address for author one
%\curraddr..... \endcurraddr
%Current address for author one
\email burdet\@cptsu2.univ-mrs.f r\endemail

%  Math Subject Classifications 
\subjclass  53 C 99 \endsubjclass

\abstract    This note describes the construction of c\ $\cup$\ p\ -invariant differential
operators on statistical manifolds, i.e. of operators canonically associated to a geometry  which
synthetizes the properties of conformal and projective geometries.    \endabstract

%  \thanks will become a 1st page footnote.
%  Use \endgraf to indicate a new paragraph; a blank line or \par will
%  be recognized as an error.
%  Don't type a period at the end; it will be supplied.

% \thanks The first author was supported in part by NSF
% Grant \#000000.\endgraf
% The final versi on of this paper will be submitted for
% publication elsewhere \endthanks

\endtopmatter

\document
\head 1. Introduction \endhead

\par Let $(S,g,t)$ be  a statistical manifold , i.e. a Riemannian manifold $(S,g)$
endowed with a fully symmetric 3 covariant tensor field $t$, the $skewness \  \ 
tensor$. This  tensor has been introduced to formalize the notion of statistical
curvature through the introduction of  a family $\lbrace \dca \rbrace_{\alpha \in \bb
R}$ of affine connections defined as follows  \cite {1}\cite {2}  
 $$ {\dca}_\roman x Y =  {\dczer} _\roman x Y - {\alpha\over 2}  (t\cdot g^{-1}) (X,Y) ,
\eqno(1.1)$$  for any couple $X, Y$ of vector fields over $S$, $\dczer{}$ denoting the
canonical Levi-Civita  connection and the dot (.) an obvious contraction. In fact
$\dca{}$ is the unique torsionfree connection satisfying
$$\dca \, g = \, \alpha \, t \, , \eqno(1.2)$$
and $(\dca , g)$ form a Codazzi pair  \cite {3}, i.e.  Codazzi equations are satisfied
$$ ({\dca}_ \roman x g)(Y,Z) = ({\dca}_\roman y g)(X,Z) . $$
In a previous paper  \cite {4}, we have dealt with the so-called  {\it c$\cup$p-geometry},
where the c$\cup$p symbol  refers to a conformal implementation of the projective
geometry in which the properties of both conformal and  projective geometries are preserved,
as it can be seen from the following
  \proclaim {Definition} 
\quad Let $\alpha$ be  a real number. Statistical manifolds $(S,g,t)$  and  $(S,\tilde g,\tilde t)$
are
$\alpha$-c$\cup$p-{\it related} if, for some  positive smooth function $\eta$ over $S$,

\centerline {$\tilde g = \eta g$,\quad 
$\tilde t = \eta ( t + sym(g\otimes  \psi))$}
\medskip
\noindent  where $\psi$ is an exact  1-form satisfying  $ d(Log \ \eta)= -\alpha \, \psi.$
 \endproclaim
\bigskip   and
 \proclaim { Proposition} 
\quad For any  real $\alpha$, the families of
$\alpha$-connections of two $\alpha$-c$\cup$p-related statistical manifolds are projectively
 related, i.e. the corresponding 1-form connections satisfy the following relation 
$${\tilde \varphi}^\alpha - \varphi^\alpha =  -\alpha \,( \theta\otimes \psi +
(\psi \cdot \theta) \, Id), \eqno(1.3) $$
and the curvature 2-forms are such that 
$${\tilde\Phi}^{\alpha} - \Phi^\alpha = \alpha \, \theta \wedge  \dca \psi 
+ \alpha^2\, ( \theta \otimes \psi)\wedge  (\psi \cdot \theta). \eqno(1.4)$$
\endproclaim
\bigskip
The invariants of the c$\cup$p-structure themself are described in \cite {4}, here we want to
study the closely related question of invariance of auxilliary objects associated to the
c$\cup$p-structure.

\head 2. c$\cup$p-invariant second order differential operators \endhead 
\par As the conformal one's, the c$\cup$p-geometry depends on the (conformal)
$\eta$ factor only, so we take up the step traced by the  conformal experience by saying that  a
differential operator $  D$, depending on a choice of the couple $(g, t)$, is
{\it  c $\cup$p-invariant of type (r;s) } if the operator $\tilde D$, corresponding to the
rescaled metric $\tilde g$ and to the modified tensor  $\tilde t$, is such that  $\tilde D
\tilde f = \eta^s  D f$, where $f$ is a {\it density of weight r}, section of a suitable line
bundle, and $\tilde f = \eta^r f$, $r,s \in \bb R$.

 Let us consider the  $\alpha-$Hessian $\dca  \otimes\ d$ acting on functions on $S$, we look for
the existence of a modification of it, form invariant for the  c$\cup$p-geometry. Knowing that
 conformal invariant operators may be expressed in terms of the Levi-Civita connection
and the Ricci curvature of a metric in the conformal class \cite 5, we study the operator
 $\dca  \otimes\ d + k \  Ric^\alpha, k \in \bb R$, where  $ Ric^\alpha$ denotes the Ricci tensor
of the $\alpha$-connection.

Starting from Rel(1.4) it can be shown that the Ricci tensors of two $\alpha$-c$\cup$p-related 
connections are such that
$$\tilde Ric^\alpha = Ric^\alpha + (n-1)\alpha \lbrace \dca \psi +  \alpha  ( \psi \otimes \psi)
\rbrace, \eqno(2.1)$$
where $n = dim S$.

Then, by using Rel(1.3) one gets
$$(\tilde \nabla^{\alpha}  \otimes\ d)( \eta^r f) = \eta^r ( \nabla^{\alpha}  \otimes\ d) f  + f (
\nabla^{\alpha}  \otimes\ d)  \eta^r + d(\eta^r)\otimes\nabla^{\alpha} f + df
\otimes\nabla^{\alpha}(\eta^r), \eqno(2.2)$$ 
taking into account for $ d(Log \eta) = -\alpha\psi$, one establishes the following
 \proclaim {Proposition}  The modified Hessian  $\dca \otimes\ d + k \  Ric^\alpha$ is
$\alpha-c\cup p$-invariant  of type  (1;1) for $k = {1 \over n-1}$:

$$(\tilde \nabla^{\alpha}  \otimes\ d + {1\over {n-1}} \ \tilde Ric^\alpha) \tilde f = \eta (\dca 
\otimes\ d + {1\over {n-1}} \  Ric^\alpha) f . \eqno(2.3)$$ 
 \endproclaim  

\n {\bf Remark}   The integrability conditions of the differential system 
$(\dca  \otimes\ d + k \  Ric^\alpha)f = 0 $  are given by$$Riem^\alpha = k(Id \otimes Ric^\alpha -
Ric^\alpha \otimes Id).$$
Except for the flat case $Riem^\alpha = 0$, this property is verified to hold, by taking the trace,
for the  c$\cup$p-invariant  Hessian only, and for statistical manifolds such that
$$Riem^\alpha = {1\over {n-1}}(Id \otimes Ric^\alpha - Ric^\alpha \otimes Id).\eqno(2.4)$$
For instance this property is satisfied for the statistical manifold associated to the gaussian  law
and for the statistical manifold  of  multinomial  distributions.
\bigskip

\head 3. The  c$\cup$p-Laplacian \endhead
\bigskip \bigskip 
\par Let us consider the trace of the Ricci modified Hessian operator
$$  g^{-1}(\dca \otimes\ d  + k \  Ric^\alpha)  \ =  g^{-1}(\dca  \otimes\ d) + k \ R^\alpha $$
where $R^\alpha $ denotes the {\it scalar  $\alpha$-curvature} $R^\alpha = g^{-1}(Ric^\alpha)$.
This expression can be written by using the $\alpha-Laplacian$ = div$^\alpha$   grad, denoted by 
$\tra$, but we have to take into account for the non-metricity of the connection; so we are led to
give the
\proclaim { Definition } 
\quad The second order differential operator given by
$$\triangle_{c\cup p} \ =  g^{-1}(\dca  \otimes\ d + {1\over {n-1}} \ Ric^\alpha ) =  \tra + \alpha
g^{-1}. t .g^{-1}. d + {1\over {n-1}} \ R^\alpha , \eqno(3.1)$$ 
where each  point denotes an obvious
contraction due to the  symmetry of implied tensor fields, is called the $ c \cup p\
- Laplacian.$ \endproclaim
Indeed this  operator is  associated to a c$\cup$p\ -structure canonically as it is clear from the
following

\proclaim { Proposition } 
$$\tilde\triangle_{c\cup p} \ \tilde f =  \triangle_{c\cup p} \  f , \  \  (\tilde f = \eta f), 
\eqno(3.2)$$ 
$the  c\cup p-Laplacian \ is \ a \ c\cup p-invariant \ differential \ operator \ of \ type \ (1;0)
$.
   \endproclaim

This property stems from Rel(2.3) by noting that ${\tilde g}^{-1} =
{1\over\eta} g^{-1}.$

\head 4. Quasi-c$\cup$p-invariant non-linear operators \endhead

\par We now wish to construct a non-linearly modified c$\cup$p- Laplacien
through a self-interaction term $f^a$ with a coupling  potential $\lambda$, function on $S$. We can
use property (3.2) to get immediately the

\proclaim { Proposition } $$\tilde\triangle_{c\cup p} \ \tilde f  + \tilde \lambda \ \tilde f ^a
= \triangle_{c\cup p} \  f   +  \lambda \ f ^a, \  \  (\tilde f = \eta f , \tilde \lambda = \lambda
/ \eta ^a , a = constante \not = 0).  \eqno(4.1)$$
 \endproclaim
It is worthnoticing that there is no constraint over the non-linearity coefficient $a$
contrarily to the conformal case for which  $a$ depends on the dimension of the manifold \cite 6.

%\head 5. Conclusion \endhead

\enddocument

\item{\hbox to\parindent{\enskip }\hfill}S-i. Amari,
O.E. Barndorff-Nielsen, R.E. Kass,
S.L. Lauritzen, C.R. Rao, {\sl Differential Geometry  in Statistical
Inference},  IMS Lecture No\-tes-Monographs Series,{\bf 10}, S.S.
Gupta Ed. (1987).
\medskip

\enddocument